\documentclass[preprint,sort&compress,12pt]{elsarticle}
\usepackage{amsmath}
\usepackage{mathrsfs}
\usepackage{stmaryrd}
\usepackage{mathrsfs}
\usepackage{bm}
\usepackage{amsmath}
\usepackage{amsfonts}
\usepackage{amssymb}
\usepackage{amsthm}
\journal{\quad}

\newcommand{\mf}{\mathbf}
\newcommand{\mm}{\mathrm}

\begin{document}
\begin{frontmatter}
\title{An Improved Result on Rayleigh--Taylor Instability of \\
Nonhomogeneous Incompressible Viscous Flows
}

%%%%%
\author[FJ]{Fei Jiang%\corref{cor1}
}
\ead{jiangfei0591@163.com}
%\author[sJ]{Song Jiang}
%\ead{jiang@iapcm.ac.cn}
\address[FJ]{College of Mathematics and Computer Science, Fuzhou University, Fuzhou, 350108, China.}
%\address[sJ]{Institute of Applied Physics and Computational Mathematics, %P.O. Box 8009,
% Beijing, 100088, China.}
\begin{abstract}
In [F. Jiang, S. Jiang, On instability and stability of three-dimensional gravity driven viscous flows in a bounded domain, Adv. Math., 264 (2014) 831--863], Jiang et.al.
investigated the instability of Rayleigh--Taylor steady-state of a three-dimensional nonhomogeneous incompressible viscous flow driven by gravity in a bounded domain $\Omega$ of class $C^2$. In particular, they proved the steady-state is nonlinearly unstable under a restrictive condition of that the derivative function of steady density possesses  a positive lower bound.
In this article, by exploiting a standard energy functional and more-refined analysis of  error estimates in the bootstrap argument,
we can show the nonlinear instability result  without the restrictive condition.
\end{abstract}

\begin{keyword}
 Navier--Stokes equations, steady state solutions, Rayleigh--Taylor instability.
 %\MSC[2000] 76D05 \sep  76E09.
%(2000 is the default)

\end{keyword}
\end{frontmatter}

% \linenumbers

\newtheorem{thm}{Theorem}[section]
\newtheorem{lem}{Lemma}[section]
\newtheorem{pro}{Proposition}[section]
\newtheorem{cor}{Corollary}[section]
\newproof{pf}{Proof}
\newdefinition{rem}{Remark}[section]
\newtheorem{definition}{Definition}[section]

\section{Introduction}
\label{Intro} \numberwithin{equation}{section}

The motion of a three-dimensional (3D) nonhomogeneous incompressible
viscous fluid in the presence of a uniform gravitational field in a
bounded domain $\Omega\subset {\mathbb R}^3$ of $C^2$-class is governed by the following
Navier--Stokes equations:
\begin{equation}\label{0101}\left\{\begin{array}{l}
 \rho_t+{  v}\cdot\nabla \rho=0,\\[1mm]
\rho {v}_t+\rho {  v}\cdot\nabla { v}+\nabla p=\mu\Delta {v}-g\rho {e}_3,\\[1mm]
\mathrm{div} {v}=0,\end{array}\right.\end{equation}
where the unknowns $\rho:=\rho(t, {x})$, ${v}:={v}(t, {x})$ and $p:=p(t, {x})$ denote the density,
velocity and pressure of the fluid, respectively;
$\mu>0$ stands for the coefficient of shear viscosity, $g>0$ for the gravitational
constant, $ {e}_3=(0,0,1)$ for the vertical unit vector, and $-g{e}_3$ for the gravitational force.
In the system \eqref{0101} the equation \eqref{0101}$_1$ is the continuity equation, while
\eqref{0101}$_2$ describes the balance law of momentum.

We studied  the instability  of the following  Rayleigh--Taylor (RT) steady-state to the system \eqref{0101} as in \cite{JFJSOOADvace}:
\begin{equation}\label{steadystate}
 {v}(t,{x})\equiv  {0}\mbox{ and } \nabla \bar{p}=-g\bar{\rho} {e}_3\quad \mbox{ in }\Omega,
\end{equation}where
\begin{equation}\label{0102}
\bar{\rho}\in C^{2}(\bar{\Omega}),\ \inf_{ {x}\in
 \Omega}\{\bar{\rho}( {x})\}>0\mbox{ and
}\partial_{x_3}
\bar{\rho}>0\mbox{ for some }{x}_0\in \Omega.
\end{equation}
It is easy to show that the steady density $\bar{\rho}$ only depends on $x_3$, the third component of ${x}$. Hence we denote
$\bar{\rho}':=\partial_{x_3}\bar{\rho}$ for simplicity.
Moreover, we can compute out the associated steady pressure $\bar{p}$ determined by $\bar{\rho}$.
The third condition posed on $\bar{\rho}$ in \eqref{0102} means that there is a region in which the RT density profile has
larger density with increasing $x_3$ (height), thus leading to the nonlinear
RT instability as shown in Theorem \ref{thm:0101new} below. RT  instability is well known as gravity-driven instability in fluids
when a heavy fluid is on top of a light one.

To investigate the RT instability of the system \eqref{0101} around the steady-state \eqref{steadystate}, we denote the perturbation by
$$ \varrho=\rho -\bar{\rho},\quad  u= {v}- {0},\quad q=p-\bar{p},$$
then, $(\varrho , u,q)$ satisfies the perturbed equations:
\begin{equation}\label{0105}\left\{\begin{array}{l}
\varrho_t+{  u}\cdot\nabla (\varrho+\bar{\rho})=0, \\[1mm]
(\varrho+\bar{\rho}){ u}_t+(\varrho+\bar{\rho}){ u}\cdot\nabla
{ u}+\nabla q=\mu\Delta{ {{u}}}-g\varrho {e}_3,\\[1mm]
 \mathrm{div} u=0.\end{array}\right.  \end{equation}
To complete the statement of the perturbed problem, we specify the
initial and boundary conditions:
\begin{equation}\label{0106}
(\varrho,{  u})|_{t=0}=(\varrho_0,{  u}_0)\quad\mbox{in } \Omega
\end{equation}
and
\begin{equation}\label{0107}
{  u}|_{\partial\Omega}= 0\quad \mbox{ for any }t>0.
\end{equation}
Moreover, the initial data should satisfy the compatibility
conditions $u_0|_{\partial\Omega}=0$ and  $\mathrm{div} u_0=0$.
If we linearize the equations \eqref{0105} around the steady-state
$(\bar{\rho}, {0})$, then the resulting linearized equations read as
\begin{equation}\label{0108}
\left\{\begin{array}{ll}
 \varrho_t+ \bar{\rho}'u_3=0, \\[1mm]
  \bar{\rho} u_t +\nabla q=\mu\Delta{ {{u}}}-g\varrho{e}_3,\\[1mm]
  \mathrm{div} u=0,
\end{array}\right.\end{equation}
where $u_3$ denotes the third component of $ u$.

 Here we briefly introduce the research progress for RT instability of
continuous flows, please refer to \cite{PJSGOI5,GYTI1,wang2011viscous,GYTI2} for incompressible and compressible stratified fluids,
and \cite{KMSMSP,HRWP,JFJSWWWOA,DRJFJS} for stratified MHD fluids. Instability of the linearized problem (i.e. linear instability) for an incompressible fluid was first introduced by Rayleigh in 1883
\cite{RLAP}.
In 2003, Hwang and Guo \cite{HHJGY} proved the nonlinear RT
 instability of $\|(\varrho, u)\|_{L^2(\Omega)}$ in the sense of Hadamard
 for a 2D nonhomogeneous incompressible inviscid fluid (i.e. $\mu=0$ in the equations \eqref{0105})
 with boundary condition $ u\cdot  {n}|_{\partial\Omega}=0$, where
$\Omega=\{(x_1,x_2)\in \mathbb{R}^2~|~-l<x_2<m\}$ and $ {n}$ denotes the
 outer normal vector to $\partial\Omega$. Jiang et.al. \cite{NJTSC2} showed the nonlinear RT
 instability of $\| {u}_3\|_{L^2(\mathbb{R}^3)}$ for the Cauchy problem of \eqref{0105}
 in the sense of Lipschitz structure, and further gave the nonlinear RT
 instability of $\|u_3\|_{L^2(\Omega)}$ in \cite{JFJSWWWN} in the sense of Hadamard in a
unbounded horizontal period domain $\Omega$.

 Recently, for a general bounded domain $\Omega$,
Jiang et.al.
 showed that the steady-state \eqref{steadystate} to the linearized problem \eqref{0105}--\eqref{0107} is linear unstable (i.e., the linear solution grows in
time in $H^2(\Omega)$) by constructing a standard energy functional for  the time-independent system of
\eqref{0108}
and exploiting a modified variational
method.
Based on the linear instability result, they further showed the
nonlinear instability of the  perturbed problem \eqref{0105}--\eqref{0107} by a bootstrap technique under the following restrictive condition
(i.e., the derivative function of steady-density enjoys a positive lower bound):
 \begin{eqnarray}\label{nnnn0103}\inf_{ {x}\in \Omega}\{\bar{\rho}'( {x})\} >0.
\end{eqnarray}
The bootstrap technique has its origins in the
 paper of  Guo and Strauss \cite{GYSWIC,GYSWICNonlinea}. It was developed by
 Friedlander et.al. \cite{FSSWVMNA}, and widely quoted in the nonlinear instability literature, see \cite{GYHYQCQ,FSNPVVNC,GRnierEmanu,FrrVishikM,BIIS413,GDOS2,LinzhuwuSS} for examples.
However the Duhamel's principle in the standard bootstrap argument can not be directly applied to show
the nonlinear instability of the problem \eqref{0105}--\eqref{0107}, see \cite{JFJSOOADvace} for the details.
 To circumvent this obstacle,
Jiang et.al. used some specific energy error estimates to replace Duhamel's principle, in which the key step is to deduce an  error estimate for $(\varrho^{\mm{d}},u^\mm{d})$ in ${L^2(\Omega)}$ (i.e the $L^2(\Omega)$-norm of difference between a nonlinear
solution $(\varrho^\delta,u^\delta)$ to the problem \eqref{0105}--\eqref{0107} and
a linear solution $(\varrho^{\mm{a}},u^\mm{a})$ to the problem \eqref{0106}--\eqref{0108})
in the bootstrap technique. To this purpose,
they introduced a new energy functional under the condition \eqref{nnnn0103}
to avoid the integrand term
 \begin{equation}\label{prtialttwotimes}
\int_0^t<((\varrho^\delta+\bar{\rho}){u}_\tau^{\mm{d}})_\tau,{u}_\tau^{\mm{d}}>\mm{d}\tau,
 \end{equation}
 since the energy estimate of Gronwall-type (see \eqref{0204}) does not directly offer any estimate
 for the term $((\varrho^\delta+\bar{\rho}){u}_\tau^{\mm{d}})_\tau $. Here $<\cdot,\cdot>$ denotes the corresponding dual product between the two spaces $H^{-1}_\sigma(\Omega)$ and $H_\sigma^1(\Omega)$, and $H^{-1}_\sigma(\Omega)$ represents the dual space of
$H_\sigma^1(\Omega):=\{u\in H^1_0(\Omega)~|~\mm{div}u=0\}$.
Using the new energy functional, they can get a sharp growth rate $\Lambda$ of any linear solution $(\varrho,{u})$ in the
 norm ``$\sqrt{\|\varrho\|_{L^2 }^2+\|{u}\|_{L^2(\Omega)}^2}$''.
Thus, applying this property to the process of specific energy error estimates, they easily obtained the desired error estimate, and thus showed the nonlinear instability.

This article is devoted to canceling the condition \eqref{nnnn0103} in the proof of nonlinear instability
in \cite{JFJSOOADvace}. More precisely, we establish the following improved  result by using a standard energy functional
and more-refined analysis techniques to deduce the  error estimate for $\|(\varrho^{\mm{d}},u^\mm{d})\|_{L^2(\Omega)}$
in the bootstrap argument, which will be showed in Section \ref{sec:03}.
\begin{thm}\label{thm:0101new}
Assume that the steady density $\bar{\rho}$ satisfies \eqref{0102}.
 Then, the steady-state \eqref{steadystate} of the system \eqref{0105}--\eqref{0107}
 is unstable in the Hadamard sense, that is, there are positive constants $\Lambda$, $m_0$, $\varepsilon$ and $\delta_0$,
 and functions $(\bar{\varrho}_0,\bar{{u}}_0)\in H^2(\Omega) \times H^2(\Omega) $,
such that for any $\delta\in (0,\delta_0)$ and initial data
 $(\varrho_0,{u}_0):=(\delta\bar{\varrho}_0,\delta\bar{{u}}_0)$
there is a unique strong solution $({\varrho},{u})\in C^0([0,{T^{\max}}),H^2(\Omega) \times H^2(\Omega) )$ of \eqref{0105}--\eqref{0107}
with a associated pressure $q\in C^0([0,{T^{\max}}),H^1(\Omega)) $, such that
\begin{equation*}\label{0115}\|\varrho(T^\delta)\|_{L^2(\Omega)},\ \|(u_1,{u}_2)(T^\delta)\|_{L^2(\Omega)},\ \|{u}_3(T^\delta)\|_{L^2(\Omega) }\geq {\varepsilon}\;
\end{equation*}
for some escape time $T^\delta:=\frac{1}{\Lambda}\mm{ln}\frac{2\varepsilon}{m_0\delta}\in
(0,{T^{\max}})$, where $T^{\max}$ denotes the maximal time of existence of the solution
$(\varrho,{u})$.
\end{thm}
By virtue of \cite{JFJSOOADvace}, the key step in the proof of Theorem \ref{thm:0101new} is to establish a  error estimate
\begin{equation}\label{nesweoesimt}\|(\varrho^{\mm{d}},u^\mm{d})\|_{L^2(\Omega)}\leq C\delta^3e^{3\Lambda t}\mbox{ for some constant }C
 \end{equation}without the restrictive condition \eqref{nnnn0103} (i.e., Lemma \ref{erroestimate}).
Here we sketch the main idea in the proof of \eqref{nesweoesimt} without \eqref{nnnn0103}. In view of the property of standard energy functional (see \eqref{0111nn}),   $\Lambda$
is also  a sharp growth rate of any linear solution $(\varrho,{u})$ in the
 norm ``$\sqrt{\|\varrho\|_{L^2 }^2+\|{u}\|_{H^2(\Omega)}^2}$'', see \cite[Proposition 3.3.]{JFJSOOADvace}.  When applying the sharp growth rate of the standard energy functional  to the  process of specific energy error estimates, we need to deal with the difficulty arising from the term \eqref{prtialttwotimes}.
 However, by a classical regularization method, we can show  that
\begin{equation*}\begin{aligned}\label{}
&2\int_0^t<((\varrho^\delta+\bar{\rho}){u}_\tau^{\mm{d}})_\tau,{u}_\tau^{\mm{d}}>\mm{d}\tau
\\
&=\int_\Omega  (\varrho^\delta+\bar{\rho}) | u_{t}^{\mm{d}}(t)|^2\mm{d}x-\int_\Omega ( \varrho^\delta(0)+\bar{\rho}) | u_{t}^{\mm{d}}(0)|^2\mm{d}x +
\int_0^t\int_\Omega \varrho_\tau^\delta|{u}_\tau^{\mm{d}}|^2
\mm{d}{x}\mm{d}\tau,
\end{aligned}
\end{equation*}
Then, we can deduce from the error equations (see \eqref{h0407})
that \begin{equation*}
\begin{aligned}
&\|\sqrt{\varrho^\delta+\bar{\rho}} u_t^\mm{d}(t)\|^2_{L^2 }+2\mu\int_0^t\|\nabla   u_\tau ^\mm{d}\|^2_{L^2 }\mm{d}\tau =\int g\bar{\rho}'|{u}^{\mathrm{d}}_3(t)|^2\mm{d} x
 +R_1+R_2(t),
\end{aligned}\end{equation*}
where  the two higher-order terms $R_1$ and $R_2(t)$ (see \eqref{defineitoofR1} and \eqref{defineitoofR2} for their definitions) can be controlled by $\delta^3 e^{3\Lambda t}$. Using the definition of  sharp growth rate, we can further estimate that
\begin{equation*}\label{}
\begin{aligned}
&\|\sqrt{\varrho^\delta+\bar{\rho}} u_t^\mm{d}(t)\|^2_{L^2}+2\mu\int_0^t\|\nabla   u_\tau ^\mm{d}\|^2_{L^2}\mm{d}\tau  \\
&\leq {\Lambda^2} \|\sqrt{\varrho^\delta+\bar{\rho}} u^\mm{d}(t)\|_{L^2}
+ {\Lambda}\mu\|\nabla u^\mm{d}(t)\|^2_{L^2}+C\delta^3 e^{3\Lambda t}.
\end{aligned}\end{equation*}Based on the estimate above, by more-refined analysis, we can further obtain the following  Gronwall's inequality
\begin{equation*}
\begin{aligned}
& \frac{\mm{d}}{\mm{d}t} \|\sqrt{\varrho^\delta+\bar{\rho}} u^\mm{d}(t)\|^2_{L^2}+ \mu\| \nabla u^\mm{d}(t)\|_{L^2}^2\\
&\leq 2\Lambda\left( \|\sqrt{\varrho^
\delta+\bar{\rho}} u^\mm{d}(t)\|^2_{L^2}
 + \mu\int_0^t\| \nabla u^\mm{d} \|_{L^2}^2
\mathrm{d}\tau\right)
+C\delta^3 e^{3\Lambda t}.
\end{aligned}
\end{equation*}
Since  $\varrho^\delta+\bar{\rho}$  possesses  a positive lower bound, thus we
immediately get the desired error estimate \eqref{nesweoesimt} from the Gronwall's inequality above and the mass equation. We mention that  Jiang et.al.  \cite{JFJSOOADvace} used another energy functional and the restrictive condition \eqref{nnnn0103} to deduce
 the following  Gronwall's inequality
   \begin{equation*}
\begin{aligned}
\frac{\mm{d}}{\mm{d}t}\int\left(\frac{|\varrho^{\mathrm{d}}|^2}{\bar{\rho}'} + \frac{\bar{\rho}| {u}^{\mathrm{d}}|^2}{g}\right)\mm{d}\mf{x}
\leq 2\Lambda\int\left(\frac{|\varrho^{\mathrm{d}}|^2}{\bar{\rho}'} + \frac{\bar{\rho}| {u}^{\mathrm{d}}|^2}{g}\right)\mm{d}\mf{x}
+C\delta^3e^{3\Lambda t}    \end{aligned}  \end{equation*}
and thus got \eqref{nesweoesimt} under \eqref{nnnn0103}.

Finally, we end this section by explaining the notations
used throughout the rest of this article. For simplicity, we
drop the domain $\Omega$ in Sobolve spaces and the corresponding norms as well as in
integrands over $\Omega$, for example,
\begin{equation*}  \begin{aligned}&
L^p:=L^p(\Omega),\
{H}^k:=W^{k,2}(\Omega),\ H^1_\sigma:=H^1_\sigma(\Omega),\ \int:=\int_\Omega .  \end{aligned}\end{equation*}In addition, we denote $I_T:=(0,T)$ and $\bar{I}_{T}:=[0,T]$ for simplicity.
%% Similarly, we also drop the domain $\Omega$ in the notation of integrand over $\Omega$, i.e.,

\section{Preliminaries}\label{sec:05}
This section is devoted to introduction of two auxiliary results,
which were established in \cite{JFJSOOADvace} and will be used to prove
Theorem \ref{thm:0101new} in next section.
The first result is about the  instability result of the linearized problem  \eqref{0106}--\eqref{0108}.
\begin{pro}\label{thm:0101}
 Assume that the steady density $\bar{\rho}$ satisfies \eqref{0102}.
Then the steady-state \eqref{steadystate} of the linearized system \eqref{0106}--\eqref{0108}
 is unstable. That is, there exists an unstable solution
$$(\mf{\varrho},{u},{q}):=e^{\Lambda t}(-\bar{\rho}'\tilde{v}_3/\Lambda,\tilde{{v}},\tilde{p})$$
 to \eqref{0106}--\eqref{0108}, where $(\tilde{{v}},\tilde{p})\in H^2 \times H^1 $
 solves the following boundary problem
 \begin{equation*}
\left\{
                              \begin{array}{ll}
\Lambda^2\bar{\rho}\tilde{{v}}+\Lambda\nabla \tilde{{p}}
=\Lambda\mu\Delta\tilde{{{v}}}+g
\bar{\rho}'\tilde{v}_3 {e}_3,\\[1mm]
\mathrm{div}\, \tilde{{v}}=0,\quad
\tilde{{v}}|_{\partial\Omega}=0
\end{array}
                            \right.
\end{equation*} with the positive constant growth rate $\Lambda$ defined by
\begin{equation}\begin{aligned}
\label{0111nn} \Lambda^2=\sup_{\tilde{{w}}\in
{H}^\sigma_0 }\frac{g \int  \bar{\rho}'\tilde{{w}}_3^2\mm{d}{x}-\Lambda\mu
\int_\Omega|\nabla
\tilde{{w}}|^2\mm{d}{x}}{\int \bar{\rho}|\tilde{{w}}|^2\mathrm{d}{x}}
.
\end{aligned}\end{equation}
 Moreover, $\tilde{{v}}$ satisfies
 $\tilde{{v}}_3\equiv\!\!\!\!\!\!/\ 0$, $
\tilde{{v}}_1^2+\tilde{{v}}_2^2\equiv\!\!\!\!\!\!/\ 0$ and
\begin{equation}\label{unsatisfiability}
\bar{\rho}'\tilde{v}_3\equiv\!\!\!\!\!\!/\ 0,
\end{equation}
where $\tilde{v}_i$ denotes the $i$-th component of $\tilde{v}$ for $i=1$, $2$, $3$.
\end{pro}
\begin{rem}The linear instability was showed in \cite[Theorem 1.1]{JFJSOOADvace}
except
\eqref{unsatisfiability}. However, we can easily get \eqref{unsatisfiability} by contradiction. Suppose
that  $\bar{\rho}'\tilde{v}_3\equiv 0$, then
\begin{equation*}\begin{aligned}\label{}
0<\Lambda^2 =&\frac{ g\int\bar{\rho}'\tilde{{v}}_{3}^2\mm{d}  x-\Lambda\mu\int|\nabla
\tilde{ v}|^2\mm{d} x}{\int \bar{\rho}|\tilde{v}|^2\mm{d}x}=-\Lambda\mu\frac{\int|\nabla
\tilde{ v}|^2\mm{d} x}{\int \bar{\rho}|\tilde{v}|^2\mm{d}x}<0 ,
\end{aligned}\end{equation*}
which contradicts. Therefore, \eqref{unsatisfiability} holds.
\end{rem}

The second result is about a local existence result of a unique strong solution to the  perturbed
problem \eqref{0105}--\eqref{0107}, which enjoys an energy estimate of Gronwall-type, see \cite[Proposition 3.3]{JFJSOOADvace} for the detailed proof.
\begin{pro} \label{pro:0401} Assume that the steady density $\bar{\rho}$ satisfies \eqref{0102}.
For any given initial data $(\varrho_0, u_0)\in H^2\times H^2$ satisfying
$\inf_{ x\in\Omega}\{\rho_0( x)\}>0$, and  the compatibility conditions
$ u_0|_{\partial\Omega}= {0}$ and $\mm{div} u_0=0$, there exist a unique strong solution
$(\varrho, u)\in C^0([0,T^{\max}),H^2\times H^2)$ to the perturbed
problem \eqref{0105}--\eqref{0107} with a associated pressure $q\in C^0([0,{T^{\max}}),H^1)$, where
  $T^{\mm{max}}$ denotes the maximal time of existence. Moreover,
\begin{enumerate}
  \item[(1)] $u_t\in C^0([0,{T^{\max}}),L^2)$ and $$0<\inf_{ x\in\Omega}\{\varrho_0( x)+\bar{\rho}\}\leq
      \inf_{ x\in\Omega}\{\varrho(t,x)+\bar{\rho}\}\leq
\sup_{ x\in\Omega}\{\varrho(t,x)+\bar{\rho}\}\leq      \sup_{ x\in\Omega}\{\varrho_0( x)+\bar{\rho}\}<+\infty $$ for any $t\in [0,T^{\max})$.
  \item[(2)]there is a constant $\bar{\delta}_0\in (0,1)$, such that
if $\mathcal{E}(t)\leq \bar{\delta}_0$ on some interval $\bar{I}_T\subset [0,T^{\max})$, then the strong solution satisfies
\begin{equation}\label{0204}
\begin{aligned}&
\mathcal{E}^2(t)+\|( u_{t},\nabla q)(t)\|_{L^2 }^2 +\int_0^t\|(\nabla u, u_\tau,\nabla
 u_\tau)\|_{L^2 }^2\mm{d}\tau\\
  &\leq C_1\left(\mathcal{E}_0^2+\int_0^t\|(\varrho, u)\|_{L^2 }^2\mm{d}\tau\right)
 \end{aligned}
\end{equation}
for any $t\in \bar{I}_T$,
where
we have defined that
\begin{equation*}
\begin{aligned}
&\mathcal{E}(t):={\mathcal{E}}((\varrho, u)(t))
=\sqrt{
\|\varrho(t)\|_{L^2 }^2+\| u(t)\|_{H^2 }^2},\\
&\mathcal{E}_0:={\mathcal{E}}((\varrho, u)(0))
=\sqrt{
\|\varrho_0\|_{L^2 }^2+\| u_0\|_{H^2 }^2},
\end{aligned}
 \end{equation*}
 and the constant $C_1>0$ only depends on $\mu$, $g$, $\bar{\rho}$ and $\Omega$.
\end{enumerate}
\end{pro}

\section{Proof of Theorem \ref{thm:0101new}}\label{sec:03}
Now we are in a position to prove Theorem \ref{thm:0101new}. First,
 in view of Proposition \ref{thm:0101}, we can construct a linear solution
\begin{equation}\label{0501}
\left(\varrho^\mm{l},
{ u}^\mm{l}\right)=e^{{\Lambda t}}
\left(\bar{\varrho}_0,
\bar{ u}_0\right)\in H^2 \times  (H^2\cap H_\sigma^1)
\end{equation}
to the equations \eqref{0106}  with a associated pressure $q^\mm{l}=e^{\Lambda t}\bar{q}_0$, where $\bar{q}_0\in  H^1$, and
$(\bar{\varrho}_0,
\bar{ u}_0)\in H^2\times( H^2\cap H_\sigma^1)$ satisfy
\begin{eqnarray}\label{n0502}
&&\|\bar{\varrho}_0\|_{L^2 }\|{\bar{u}}_{03}\|_{L^2 }\|({\bar{u}}_{01},\bar{u}_{02})\|_{L^2 }>0,\\
&&\nonumber {\mathcal{E}}((\bar{\varrho}_0,\bar{{u}}_0))
=\sqrt{
\|\bar{\varrho}_0\|_{L^2 }^2+\|\bar{ u}_0\|_{H^2 }^2}=1,
\end{eqnarray}
where $\bar{u}_{0i}$  stands for the $i$-th component of $\bar{{u}}_0$ for $i=1$, $2$ and $3$.

Denote
$(\varrho_0^\delta,{u}_0^\delta)
:=\delta (\bar{\varrho}_0,\bar{{u}}_0)$, and $C_2:=\|(\bar{\varrho}_0,\bar{{u}}_0 )\|_{L^2}$.
Keeping in mind that the condition $\inf_{ x\in\Omega}\{\bar{\rho}( x)\}>0$ and the embedding $H^2\hookrightarrow L^\infty$, we can choose a sufficiently small $\tilde{\delta}\in (0,1)$, such that
\begin{equation*}\frac{\inf_{ x\in\Omega}\{\bar{\rho}( x)\}}{2}\leq
\inf_{ x\in\Omega}\{\varrho_0^\delta( x)+\bar{\rho}( x)\} \mbox{ for any }\delta\in (0,\tilde{\delta}).
\end{equation*}
Thus, by virtue of Proposition \ref{pro:0401},   for any $\delta<\tilde{\delta}$, there exists a unique local solution
$(\varrho^\delta, u^\delta)\in C^0([0,T^{\max}),H^2\times H^2)$ to the perturbed problem \eqref{0105}--\eqref{0107}
with a associated pressure $q^\delta\in C^0([0,{T^{\max}}),H^1)$, emanating
from the initial data $(\varrho_0^\delta, {u}_0^\delta)$ with ${\mathcal{E}}((\varrho_0^\delta, {u}_0^\delta))=\delta$, where $T^{\mm{max}}$ denotes the maximal time of existence.
Moreover,
\begin{equation}\label{infdensity}
0<\frac{\inf_{ x\in\Omega}\{\bar{\rho}( x)\}}{2}\leq
\inf_{ x\in\Omega}\{\varrho^\delta (t,x)+\bar{\rho}\}
\end{equation}and
\begin{equation}\label{uppersdensty}
 \sup_{ x\in\Omega}\{\varrho^\delta(t,x)+\bar{\rho}\}\leq \sup_{ x\in\Omega}\{\bar{\varrho}_0( x)+\bar{\rho}\}\leq C_3\|\bar{\varrho}_0\|_{H^2}+\|\bar{\rho}\|_{L^\infty}
\end{equation} for any $t\in [0,T^{\max})$, where $C_3$ is the constant from the imbedding
$H^2\hookrightarrow {L^\infty}$.

Let $C_1>0$ and $\bar{\delta}_0>0$ be the same constants as in Proposition \ref{pro:0401}, and $\varepsilon_0\in (0,1)$ be a constant, which will be defined in \eqref{defined}.
Denote
$\delta_0=\min\{\tilde{\delta},\bar{\delta}_0\}$, for given $\delta\in (0,\delta_0)$, we define
 \begin{equation}\label{times}
 T^{\delta}:=\frac{1}{\Lambda}\mm{ln}\frac{2\varepsilon_0}{\delta}>0,\quad\mbox{i.e.,}\;
 \delta e^{\Lambda T^\delta}=2\varepsilon_0, \end{equation}
 \begin{equation*}
T^*:=\sup\left\{t\in I_{T^{\max}}\left|~{\mathcal{E}}((\varrho^\delta,
{ u}^\delta )(t))\leq {\delta_0}\right.\right\}>0\end{equation*}
 and
   \begin{equation*}
   T^{**}:=\sup\left\{t\in I_{T^{\max}}\left|~\left\|\left(\varrho^\delta,
{ u}^\delta\right)(t)\right\|_{{L}^2 }\leq 2\delta C_2e^{\Lambda t}\right\}>0\right.. \end{equation*}
Then $T^*$ and $T^{**}$ may be finite, and furthermore,
 \begin{eqnarray}\label{0502n1}
&&{\mathcal{E}}(\left(\varrho^\delta, { u}^\delta\right)(T^*))={\delta_0}\quad\mbox{ if }T^*<\infty ,\\
\label{0502n111}  &&  \left\|\left(\varrho^\delta, { u}^\delta\right)(T^{**})\right\|_{{L}^2 }
=2\delta C_2e^{\Lambda T^{**}}\quad\mbox{ if }T^{**}<T^{\max}.
\end{eqnarray}
Now, we denote ${T}_{\min}:= \min\{T^\delta,T^*,T^{**}\}$, then for all $t\in
\bar{I}_{ {T}_{\min}}$, we deduce from the estimate \eqref{0204} and the
definitions of $T^*$ and $T^{**}$ that
\begin{align}
&{\mathcal{E}}^2\big( (\varrho^\delta, { u}^\delta)(t)\big) +\| u_{t}^\delta(t)\|_{L^2 }^2+\int_0^t
\|\nabla u_\tau^\delta\|_{L^2 }^2\mm{d}\tau\nonumber\\
&\leq C_1 \delta^2 {\mathcal{E}}^2(\left(\bar{\varrho}_0,
\bar{ u}_0 \right))+C_1\int_0^t\left\|\left(\varrho^\delta,
{ u}^\delta \right)\right\|_{L^2 }^2\mm{d}\tau\nonumber\\
&\label{0503}\leq  C_1\delta^2+4 C_1C_2^2\delta^2e^{2\Lambda t}/(2\Lambda)
\leq C_4\delta^2e^{2\Lambda t}
   \end{align}
where $C_4:=C_1 +4 C_1C_2^2 /(2\Lambda)$  is independent of $\delta$.

Let $(\varrho^{\mathrm{d}}, { u}^{\mathrm{d}})=(\varrho^{\delta},
{ u}^{\delta})-\delta(\varrho^{\mathrm{l}}, { u}^{\mathrm{l}})$.
Noting that $(\varrho^\mm{a},{u}^{\mm{a}}):= \delta(\varrho^{\mm{l}},{u}^{\mm{l}})\in
C^0([0,+\infty),H^2\times H^2 )$
is also a linear solution to \eqref{0106}--\eqref{0108} with the initial data
$(\varrho_0^\delta, u_0^\delta)\in H^2 \times H^2 $ and with a associated pressure $q^a=\delta
q^\mm{l}\in
 C^0([0,+\infty),H^1 )$, we find that $(\varrho^{\mathrm{d}}, { u}^{\mathrm{d}})$
satisfies the following error equations:
\begin{equation}\label{h0407}\left\{\begin{array}{ll}
  \varrho_t^{\mathrm{d}}+\bar{\rho}'{u}_3^{\mathrm{d}}= -{{ u}}^{\delta}\cdot \nabla\varrho^{\delta}, \\[1mm]
  (\varrho^\delta+\bar{\rho}) u_t^{\mathrm{d}}   -\mu \Delta  u^{\mathrm{d}}+\nabla q^{\mathrm{d}}
= f^\delta -g\varrho^{\mathrm{d}} {e}_3,\\[1mm]
 \mathrm{div} u^{\mathrm{d}}={0}
 \end{array}\right.\end{equation}
 with initial and boundary conditions
  \begin{equation*} (\varrho^{\mathrm{d}}(0),
 {{u}}^{\mathrm{d}}(0))=
{ {0}},\  u^{\mathrm{d}}|_{\partial\Omega}=0£¬
   \end{equation*}
   and  compatibility
conditions
$$u^\mm{d}(0)|_{\partial\Omega}=0,\ \mathrm{div} u^\mm{d}(0)=0,$$
   where we have defined that
   $$q^{\mm{d}}:=q^\delta-q^\mm{a}\in C^0(\bar{I}_{ {T}_{\min}},H^1)\mbox{ and   } f^\delta:=-( \varrho^{\delta}+\bar{\rho}) u^{\delta}\cdot\nabla
                              u^{\delta}-\varrho^\delta u^{\mm{a}}_t.$$

Next, we shall establish an error estimate for $(\varrho^{\mm{d}},u^{\mm{d}})$ in $L^2$-norm.
\begin{lem}\label{erroestimate}  There is a constant $C_4$, such that for all $t\in \bar{I}_{{T}_{\min}}$,
\begin{equation}\label{ereroe}
\begin{aligned}
  \| (\varrho^{\mm{d}},u^{\mm{d}})(t)\|^2_{L^2 } \leq C_4\delta^3e^{3\Lambda t}.
 \end{aligned}  \end{equation}
   \end{lem}
\begin{pf}
Recalling that
$(\varrho^{\mathrm{d}}, { u}^{\mathrm{d}})=(\varrho^\delta,u^\delta)-(\varrho^\delta,u^\delta)$, in view of the regularity of $(\varrho^\delta,u^\delta)$ and $(\varrho^\mm{a},u^{\mm{a}})$,
    we can deduce from \eqref{h0407}$_2$
 that
 for a.e. $t\in
{I}_{T_{\min}}$,
\begin{equation}\label{n0310}
\begin{aligned}
&
\frac{\mm{d}}{\mm{d}t}\int (\varrho^\delta+\bar{\rho})| u_t^{\mathrm{d}}|^2
\mathrm{d} x=2<((\varrho^\delta+\bar{\rho}){u}_t^{\mm{d}})_t,{u}_t^{\mm{d}}>-\int \varrho^\delta_t|u_t^{\mm{d}}|^2\mm{d}x\\
&=2\int (f_t^\delta-g\varrho^{\mm{d}}_t e_3) u^{\mathrm{d}}_t\mm{d}x-2\mu\int
|\nabla u^{\mathrm{d}}_t|^2
\mm{d} x-\int \varrho^\delta_t|u_t^{\mm{d}}|^2\mm{d}x,
\end{aligned}\end{equation}
and $\|\sqrt{\varrho^\delta+\bar{\rho}}u_t^{\mm{d}}\|_{L^2}\in C^0(\bar{I}_{T_{\min}})$,
please refer to   \cite[Remark 6]{CYKHOM1ufda}.
Noting that
$$\frac{\mm{d}}{\mm{d}t}\int \bar{\rho}'|{u}_3^{\mathrm{d}}|^2\mm{d}x=2\int \bar{\rho}'{u}_3^{\mathrm{d}}\partial_tu_3^{\mathrm{d}}\mm{d}x,$$
thus, using \eqref{h0407}$_1$, we can rewrite the equality \eqref{n0310} as
\begin{equation}\label{nnn0314P}
\begin{aligned}
&\frac{\mm{d}}{\mm{d}t}
\int \left[(\varrho^\delta+\bar{\rho})| u_t^{\mathrm{d}}|^2
-g\bar{\rho}'|{u}_3^{\mathrm{d}}|^2\right]
\mathrm{d} x+2\mu\int
|\nabla u^{\mathrm{d}}_t|^2
\mm{d} x\\
&=\int\left (2f_t+2gu^\delta\cdot\nabla \varrho^\delta e_3- \varrho^\delta_tu_t^{\mm{d}}\right)\cdot u_t^{\mm{d}}\mm{d}x,
\end{aligned}\end{equation}
Recalling that $u^\mm{d}_3(0)=0$,
thus, integrating (\ref{nnn0314P}) in time from $0$ to $t$, we get
\begin{equation}\label{0314}
\begin{aligned}
\|\sqrt{\varrho^\delta+\bar{\rho}} u_t^\mm{d}(t)\|^2_{L^2 }+2\mu\int_0^t\|\nabla   u_\tau ^\mm{d}\|^2_{L^2 }\mm{d}\tau  =\int g\bar{\rho}'|{u}^{\mathrm{d}}_3(t)|^2\mm{d} x
 + R_1+R_2(t),
\end{aligned}\end{equation}
where
\begin{equation}\label{defineitoofR1}
R_1=\left[\int(\varrho^\delta+\bar{\rho})| u_t^{\mathrm{d}}|^2\mathrm{d} x\right]\bigg|_{t=0}
\end{equation}
and
\begin{equation}\label{defineitoofR2}
\begin{aligned}
R_2(t)=\int_0^t
\int \left(2f_\tau+2 gu^\delta\cdot\nabla \varrho^\delta e_3- \varrho^\delta_\tau u_\tau^{\mm{d}}\right)\cdot u_\tau^{\mm{d}}\mm{d}x \mm{d}\tau.
\end{aligned}\end{equation}
Next, we control the two higher-order terms $R_1$ and $R_2(t)$ above.
In what follows,
we denote by $C$ a
generic positive constant which may depend on $\mu$, $g$, $\bar{\rho}$, $\Lambda$, $\Omega$ and $(\bar{\varrho}_0,\bar{u}_0)$.
The symbol $a\lesssim b$ means that $a\leq Cb$.

Multiplying \eqref{h0407}$_2$ by $u_t^\mm{d}$ in $L^2$, we get
$$\int(\varrho^\delta+\bar{\rho})| {u}_t^{\mathrm{d}}|^2\mm{d}x
  =\int
(f^\delta-g\varrho^\mm{d}e_3+\mu\Delta u^{\mm{d}})\cdot {u}_t^{\mathrm{d}}\mm{d}x.$$
Exploiting \eqref{infdensity} and Cauchy's inequality, we get
\begin{equation}\label{appesimtsofu12}\int(\varrho^\delta+\bar{\rho})| {u}_t^{\mathrm{d}}|^2\mm{d}x
  \lesssim
\|f^\delta-g\varrho^\mm{d}e_3\|_{L^2}^2+\|\Delta u^{\mm{d}}\|_{L^2}^2.
\end{equation}
By the definition of $u_t^{\mm{a}}$, it holds that
\begin{equation}\label{appesimtsofu}
\|\partial_{t}^ju^\mm{a}\|_{H^k }=\Lambda^j \delta e^{\Lambda t}\|\bar{u}_0\|_{H^k }\mbox{ for }0\leq k,\, j\leq 2,
\end{equation}
thus, using \eqref{uppersdensty}, \eqref{0503}, H\"older's inequality and the imbedding
$H^2\hookrightarrow L^\infty$,  we  have
\begin{equation}\label{appesofu12}
\begin{aligned}
\|f^\delta-g\varrho^\mm{d}e_3\|_{L^2}^2
  \lesssim &\|\varrho^\mm{d}\|_{L^2 }^2+
  \|(\varrho^\delta+\bar{\rho})\|_{L^\infty }^2
  \|u^\delta\|_{H^2 }^4+\|\varrho^\delta\|_{L^2 }^2\|u_t^{\mm{a}}\|_{H^2 }^2
\\
 \lesssim  &\|\varrho^\mm{d}\|_{L^2 }^2+\delta^4 e^{4\Lambda t},
 \end{aligned}
 \end{equation}
Noting that $\varrho^\mm{d}(0)=0$, $\Delta u^{\mm{d}}(0)=0$ and $\delta \in (0,1)$, chaining the estimates
\eqref{appesimtsofu12} and \eqref{appesofu12} together, and taking limit for $t\rightarrow 0$, we
immediately obtain the following estimate for the first higher-order term $R_1$:
\begin{equation}\label{estimeR1}\begin{aligned}
R_1=&\lim_{t\rightarrow 0}\int(\varrho^\delta+\bar{\rho})| {u}_t^{\mathrm{d}}(t)|^2\mm{d}x
\\
\lesssim &\lim_{t\rightarrow 0}( \|\varrho^\mm{d}(t)\|_{L^2 }^2+\|\Delta u^{\mm{d}}(t)\|_{L^2}^2+\delta^4 e^{4\Lambda t}) =\delta^4\leq \delta^3.
\end{aligned}\end{equation}

Now we turn to estimate the most complicated higher-order term $R_2(t)$.
Recalling the definition of $R_2(t)$, we see that
$$\begin{aligned}R_2(t)=&-2\int_0^t\int \left[\varrho^\delta u_{\tau\tau}^{\mm{a}}+(
\varrho^{\delta}+\bar{\rho}) {u}^\delta_\tau\cdot\nabla
u^\delta +(\varrho^{\delta}+\bar{\rho}) {u}^\delta\cdot\nabla
u^\delta_\tau
\right]\cdot
 u_\tau^{\mathrm{d}}\mm{d} x\mm{d}\tau\\
 &\quad+\int_0^t\int\left[2g{{ {u}}}^{\delta}\cdot \nabla\varrho^{\delta}e_3 -\varrho^\delta_\tau \left(2 u^{\mm{a}}_\tau+ u^{\mm{d}}_\tau+2 {u}^\delta\cdot\nabla
u^\delta\right)\right]\cdot u_\tau^{\mathrm{d}}\mm{d}x\mm{d}\tau\\
 :=&R_{2,1}(t)+R_{2,2}(t).
 \end{aligned}$$
Using \eqref{uppersdensty}, \eqref{0503},  \eqref{appesimtsofu},
H\"older's inequality
and  the imbeddings $ H^2\hookrightarrow L^\infty$ and
$H^1\hookrightarrow L^4$, the integral term $R_{2,1}(t)$ can be estimated as follows:
\begin{equation}\label{esimtaes11}
\begin{aligned}
R_{2,1}(t)\lesssim &\int_0^t
(\|\varrho^\delta \|_{L^2 }\|u^\mm{a}_{\tau\tau}\|_{H^2 }+\|(\varrho^\delta+\bar{\rho})
\|_{L^\infty }\|  u^\delta\|_{H^2 }\| u_\tau^\delta\|_{H^1 }
)\|  u^\mm{d}_\tau\|_{L^2 }\mm{d}\tau
\\
\lesssim &\int_0^t\delta^2 e^{2\Lambda \tau}  (\delta e^{\Lambda \tau}+\|\nabla u_\tau^\delta\|_{L^2 })\mm{d}\tau\\
\lesssim & \delta^3 e^{3\Lambda t}+\left(\int_0^t\delta^4 e^{4\Lambda \tau}  \mm{d}\tau \right)^{\frac{1}{2}}
\left(\int_0^t  \|\nabla u_\tau^\delta\|_{L^2 }^2\mm{d}\tau\right)^{\frac{1}{2}}\lesssim
\delta^3 e^{3\Lambda t}.
\end{aligned}
\end{equation}
To estimate the second term $R_{2,2}(t)$, we use the mass equation (i.e. $\varrho_t^\delta
=-({{ {u}}}^{\delta}\cdot \nabla\varrho^{\delta}+\bar{\rho}'{u}_3^{\delta})$) and  the formula of integration by parts to rewrite $R_{2,2}(t)$ as follows:
\begin{equation*}\label{}
\begin{aligned}
R_{2,2}(t)=&\int_0^t\int\left[\left({{ {u}}}^{\delta}\cdot \nabla\varrho^{\delta}+\bar{\rho}'{u}_3^{\delta}\right)\left(2u^{\mm{a}}_\tau+u^{\mm{d}}_\tau+
2{u}^\delta\cdot\nabla
u^\delta\right)+2g{{ {u}}}^{\delta}\cdot \nabla\varrho^{\delta}e_3 \right]\cdot u_\tau^{\mathrm{d}}\mm{d}x\mm{d}\tau\\
 =&\int_0^t\int\left[\bar{\rho}'{u}_3^{\delta}\left( 2u^{\mm{a}}_\tau+ u^{\mm{d}}_\tau+2{u}^\delta\cdot\nabla
u^\delta\right) u_\tau^{\mathrm{d}}-2g\varrho^{\delta}{{ {u}}}^{\delta}\cdot \nabla \partial_\tau u_{3}^{\mathrm{d}} \right]\mm{d}x\mm{d}\tau\\
&-2\int_0^t\int\left[\varrho^{\delta}{u}^{\delta}\cdot \nabla\left(u^\delta_\tau+{u}^\delta\cdot\nabla
u^\delta \right)\cdot u_\tau^{\mathrm{d}}+\varrho^{\delta}{u}^{\delta}\cdot\nabla u_\tau^{\mathrm{d}}\cdot \left(u^{\mm{a}}_\tau+{u}^\delta\cdot\nabla
u^\delta\right)\right]\mm{d}x\mm{d}\tau\\
=& R_{2,2,1}(t)+R_{2,2,2}(t).
 \end{aligned}\end{equation*}
Similarly to \eqref{esimtaes11}, we can estimate that
 \begin{equation}\label{estimatesr123}
\begin{aligned}
R_{2,2,1}(t)\lesssim &\int_0^t[\|u_3^\delta\|_{H^2}(\|u_\tau^{\mm{a}}\|_{L^2}+\|u_\tau^{\delta}\|_{L^2}+
\|u^\delta\|_{H^2}^2)
\|u_\tau^{\mm{d}}\|_{L^2} +
\|\varrho^\delta\|_{L^2}
\|u^\delta\|_{H^2}\|\nabla\partial_\tau u_3^{\mm{d}}\|_{L^2}]\mm{d}\tau\\
\lesssim &\int_0^t\left[\delta^3 e^{3\Lambda \tau}(1+\delta e^{\Lambda \tau})+
\delta^2e^{2\Lambda \tau}\|\nabla\partial_\tau u_3^\delta\|_{L^2}\right]\mm{d}\tau\lesssim
\delta^3 e^{3\Lambda t}(1+\delta e^{\Lambda t}),
 \end{aligned}\end{equation}
 and
\begin{equation}\label{estimatesr123n} \begin{aligned}
R_{2,2,2}(t)\lesssim&
\int_0^t\|\varrho^\delta\|_{L^\infty}
 \|u^\delta\|_{H^2}(\|\nabla u_\tau^{\delta}\|_{L^2}\|  u_\tau^{\mm{d}}\|_{L^2}+ \|u^\delta\|_{H^2}^2\|u_\tau^{\mm{d}}\|_{L^2}\\
 &\qquad +\|u_\tau^{\mm{a}}\|_{L^2}\|\nabla u_\tau^{\mm{d}}\|_{L^2}
 +\|u^\delta\|_{H^2}^2\|\nabla  u^\mm{d}_\tau\|_{L^2})\mm{d}\tau\\
\lesssim &\int_0^t\left[\delta^3 e^{3\Lambda \tau}(1+\delta e^{\Lambda \tau})+
\delta^2 e^{2\Lambda \tau}\|\nabla u_\tau^\delta\|_{L^2}\right]\mm{d}\tau\lesssim
\delta^3 e^{3\Lambda t}(1+\delta e^{\Lambda t}).
 \end{aligned}\end{equation}

By the definition of $\varepsilon_0\in (0,1)$ in \eqref{times},  \begin{equation}\label{timesd}
 \delta\leq \delta e^{\Lambda t}\leq \delta e^{\Lambda T^\delta}\leq 2\mbox{ for any }t\in \bar{I}_{T_{\min}},
\end{equation}
Thus, summing up the estimates \eqref{estimeR1}--\eqref{estimatesr123n}, we get
\begin{equation}\label{estimateforhigher}
\begin{aligned}
R_1+R_2(t)= R_1+R_{2,1}(t)+R_{2,2,1}(t)+R_{2,2,2}(t)
\lesssim \delta^3 e^{3\Lambda t},
\end{aligned}\end{equation}
which, together with \eqref{0314}, yields that
\begin{equation*}\label{}
 \|\sqrt{\varrho^\delta+\bar{\rho}} u_t^\mm{d}(t)\|^2_{L^2}+2\mu\int_0^t\|\nabla   u_\tau ^\mm{d}\|^2_{L^2}\mm{d}\tau  \leq \int g\bar{\rho}'|{u}^{\mathrm{d}}_3|^2\mm{d} x +C\delta^3 e^{3\Lambda t}.
 \end{equation*}
Thanks to \eqref{0111nn}, we have
\begin{equation*}\begin{aligned}
\int g\bar{\rho}'|{u}^{\mathrm{d}}_3|^2\mm{d} x \leq &\Lambda^2{\int\bar{\rho}|{ u}^{\mathrm{d}}|^2\mathrm{d} x} +\Lambda\mu\int|\nabla {u}^{\mathrm{d}}|^2\mm{d}x\\
=& \Lambda^2{\int(\varrho^\delta+\bar{\rho})|{ u}^{\mathrm{d}}|^2\mathrm{d} x} +\Lambda\mu\int|\nabla {u}^{\mathrm{d}}|^2\mm{d}x-\Lambda^2{\int\varrho^\delta|{ u}^{\mathrm{d}}|^2\mathrm{d} x}\\
\leq &\Lambda^2{\int(\varrho^\delta+\bar{\rho})|{ u}^{\mathrm{d}}|^2\mathrm{d} x} +\Lambda\mu\int|\nabla {u}^{\mathrm{d}}|^2\mm{d}x+C\delta^3 e^{3\Lambda t}.
\end{aligned}\end{equation*}
Chaining the previous two inequalities together, we obtain
\begin{equation}\label{new0311}
\begin{aligned}
&\|\sqrt{\varrho^\delta+\bar{\rho}} u_t^\mm{d}(t)\|^2_{L^2}+2\mu\int_0^t\|\nabla   u_\tau ^\mm{d}\|^2_{L^2}\mm{d}\tau  \\
&\leq {\Lambda^2} \|\sqrt{\varrho^\delta+\bar{\rho}} u^\mm{d}(t)\|_{L^2}^2
+ {\Lambda}\mu\|\nabla u^\mm{d}(t)\|^2_{L^2}+C\delta^3 e^{3\Lambda t}.
\end{aligned}\end{equation}

Recalling that $u^\mm{d}\in C^0(\bar{I}_{T_{\min}}, H^2)$ and $\nabla u^{\mm{d}}(0)=0$, thus, using Newton-Leibniz's formula and Cauchy--Schwarz's inequality, we find that
 \begin{equation}\begin{aligned}\label{0316}
 \Lambda \mu\|\nabla u^\mm{d}(t)\|_{L^2}^2
& =  2\Lambda\mu\int_0^t\int_{\Omega}\sum_{1\leq i,j\leq 3}\partial_{x_i} u_{j}^\mm{d} \partial_{x_i} u_{j\tau}^\mm{d} \mm{d} x\mathrm{d}\tau
 \\
& \leq \Lambda^2\mu\int_0^t\|\nabla u^\mm{d} \|_{L^2}^2\mathrm{d}\tau +\mu\int_0^t\|\nabla u_\tau^\mm{d} \|_{L^2}^2
\mathrm{d}\tau,       \end{aligned}\end{equation}
where   $u_{j\tau}^\mm{d}$ denotes the $j$-th component of
 $u_{\tau}^\mm{d}$.
Putting \eqref{new0311} and \eqref{0316} together, we have \begin{equation}\label{inequalemee}\begin{aligned}
&\frac{1}{\Lambda}\|\sqrt{\varrho^\delta+\bar{\rho}}  u_t^{\mathrm{d}}(t)\|^2_{L^2 }+
{\mu}\|\nabla u^{\mathrm{d}}(t)\|_{L^2 }^2\\
& \leq   {\Lambda}\|\sqrt{\varrho^\delta+\bar{\rho}} u^{\mathrm{d}}(t)\|^2_{L^2 }+2 {\Lambda}\mu\int_0^t\|
\nabla u^{\mathrm{d}}\|_{L^2 }^2\mm{d}\tau +C\delta^3 e^{3\Lambda t}.
\end{aligned}\end{equation}
 On the other hand,
\begin{equation*}\begin{aligned}\label{}
\frac{\mm{d}}{\mm{d}t}\|\sqrt{\varrho^\delta+\bar{\rho}} u^\mm{d} \|^2_{L^2}=&2\int
(\varrho^\delta+\bar{\rho}) u^\mm{d} \cdot  u^\mm{d}_t \mm{d} x+\int
\varrho^\delta_t |u^\mm{d}|^2 \mm{d} x\\
\leq&\frac{1}{\Lambda}\|\sqrt{(\varrho^\delta+\bar{\rho})}  u_t^\mm{d} \|^2_{L^2}
+\Lambda\|\sqrt{\varrho^\delta+ \bar{\rho}} u^\mm{d} \|^2_{L^2}+\int
\varrho^\delta_t |u^\mm{d}|^2 \mm{d} x
\end{aligned}\end{equation*}
and \begin{equation*}\label{}\begin{aligned}
\int
\varrho^\delta_t |u^\mm{d}|^2 \mm{d} x=&-\int
({{ {u}}}^{\delta}\cdot \nabla\varrho^{\delta}+\bar{\rho}'{u}_3^{\delta}) |u^\mm{d}|^2 \mm{d} x\\
=&\int
(2\varrho^{\delta}{{ {u}}}^{\delta}\cdot \nabla u^\mm{d}-
\bar{\rho}'{u}_3^{\delta}u^\mm{d})\cdot u^\mm{d}\mm{d} x\\
\lesssim & \delta^3 e^{3\Lambda t}
\end{aligned}\end{equation*}
Putting  the previous three estimates together, we get the differential inequality
\begin{equation}\label{growallsinequa}
\begin{aligned}
& \frac{\mm{d}}{\mm{d}t} \|\sqrt{\varrho^\delta+\bar{\rho}} u^\mm{d}(t)\|^2_{L^2}+ \mu\| \nabla u^\mm{d}(t)\|_{L^2}^2\\
&\leq 2\Lambda\left( \|\sqrt{\varrho^
\delta+\bar{\rho}} u^\mm{d}(t)\|^2_{L^2}
 + \mu\int_0^t\| \nabla u^\mm{d} \|_{L^2}^2
\mathrm{d}\tau\right)
+C\delta^3 e^{3\Lambda t}.
\end{aligned}
\end{equation}

Recalling that $u^{\mm{d}}=0$, thus, applying Gronwall's inequality to \eqref{growallsinequa}, one obtains
 \begin{equation}\label{estimerrvelcoity}
\begin{aligned}
 \|\sqrt{\varrho^\delta+\bar{\rho}} u^{\mathrm{d}}(t)\|^2_{L^2}+
{\mu}\int_0^t\|\nabla u^{\mathrm{d}}\|^2_{L^2}
\mm{d}\tau
\leq  e^{2\Lambda t}\int_0^t (C  \delta^3 e^{3\Lambda \tau})e^{-2\Lambda\tau}\mm{d}\tau
\lesssim
 \delta^3e^{3\Lambda t}
 \end{aligned}  \end{equation}
 for all $t\leq \bar{I}_{T_{\min}}$, which, together with \eqref{uppersdensty} and \eqref{inequalemee}, yields
 that
\begin{eqnarray}\label{uestimate1n}
\| u^{\mathrm{d}}(t)\|_{H^1  }^2+\| u_t^{\mathrm{d}}(t)\|^2_{L^2  }+
\int_0^t\|\nabla u^{\mathrm{d}}\|^2_{L^2 }\mm{d}\tau \lesssim \delta^3e^{3\Lambda t}.
\end{eqnarray}
Finally, using  the estimates
\eqref{0503}, \eqref{timesd} and \eqref{uestimate1n},
we can deduce from the equations \eqref{h0407}$_1$, that
\begin{equation}\begin{aligned}\label{erroresimts}
\|\varrho^{\mathrm{d}}(t)\|_{L^2 }\leq &\int_0^t
\|\varrho^{\mathrm{d}}_\tau\|_{L^2 }\mm{d}\tau \\
\lesssim &\int_0^t(\| u^\mm{d} \|_{H^1 }^2+\|{{ u}}^{\delta}\cdot \nabla\varrho^{\delta}\|_{L^2 })\mm{d}\tau \\
\lesssim & \int_0^t(\delta^\frac{3}{2}e^{\frac{3\Lambda}{2}\tau}+
\delta^2 e^{2\Lambda\tau})\mm{d}\tau\lesssim \delta^\frac{3}{2}e^{\frac{3\Lambda}{2} t},
\end{aligned}\end{equation}
which, together with \eqref{uestimate1n}, yields \eqref{ereroe}.
This completes the proof of Lemma \ref{erroestimate}.  \hfill$\Box$
\end{pf}

 Now, we claim that
\begin{equation}\label{n0508}
T^\delta=T_{\min},
 \end{equation}
provided that small $\varepsilon_0$ is taken to be
 \begin{equation}\label{defined}
\varepsilon_0=\min\left\{\frac{{\delta_0}}
{4},\frac{C_2^2}{8C_4},\frac{m_0^2}{C_4} \right\},
 \end{equation}where we have defined that
$m_0=: \min\{
 \|\bar{\varrho}_0\|_{L^2},\|\bar{u}_{03}\|_{L^2},
\|(\bar{u}_{01},\bar{u}_{02})\|_{L^2}\}>0$ due to \eqref{n0502}.

 Indeed, if $T^*=T_{\min}$, then
$T^*<\infty$. Moreover, from \eqref{times} and \eqref{0503} we get
 \begin{equation*}
{\mathcal{E}}(\left(\varrho^\delta,
{ u}^\delta
\right)(T^*))\leq \delta e^{\Lambda T^*}
\leq \delta e^{\Lambda T^\delta}=2\varepsilon_0<{\delta_0},
 \end{equation*}
 which contradicts with \eqref{0502n1}. On the other hand, if $T^{**}<T_{\min}$, then $T^{**}<T^*\leq T^{\mm{max}}$.
 Moreover, in view of \eqref{0501}, \eqref{times} and \eqref{ereroe}, we see that
 \begin{equation*}\begin{aligned}
 \left\|\left(\varrho^\delta,
{ u}^\delta
\right)(T^{**})\right\|_{L^2 }
\leq  & \left\|\left(\varrho^\mm{a}_{\delta},
{ u}^\mm{a}_{\delta}
\right)(T^{**})\right\|_{L^2 } +\left\|\left(\varrho^{\mathrm{d}},
{ u}^{\mathrm{d}}
\right)(T^{**})\right\|_{L^2 } \\
\leq  &\delta \left\|\left(\varrho^\mm{l},
{ u}^{\mm{l}}
\right)(T^{**})\right\|_{L^2 }+\sqrt{C_4}\delta^{3/2}e^{3\Lambda T^{**}/2} \\
\leq & \delta C_2e^{\Lambda T^{**}}+\sqrt{C_4}\delta^{3/2} e^{3\Lambda T^{**}/2}
\leq \delta e^{\Lambda T^{**}}(C_2+\sqrt{2C_4\varepsilon_0})\\
<&2\delta C_2  e^{\Lambda T^{**}},
 \end{aligned} \end{equation*}
which also contradicts with \eqref{0502n111}. Therefore, \eqref{n0508} holds.

Since $T^\delta=T_{\min}$, \eqref{ereroe} holds for $t=T^\delta$. Thus,
we can use \eqref{defined} and \eqref{ereroe} with $t=T^\delta$ to deduce that
 \begin{equation*}\begin{aligned}
 \|\varrho^{\delta}(T^\delta)\|_{L^2 }\geq &
\|\varrho^{\mathrm{a}}_{\delta}(T^{\delta})\|_{L^2 }-\|\varrho^{\mm{d}}(T^{\delta})\|_{L^2 }
= \delta\|\varrho^{\mathrm{l}}(T^{\delta})\|_{L^2 }-\|\varrho^{\mm{d}}(T^{\delta})\|_{L^2 } \\
 \geq & \delta e^{\Lambda T^\delta}\|\bar{\varrho}_{0}\|_{L^2 } -\sqrt{C_4}\delta^{3/2}e^{3\Lambda^* T^{\delta}/2}
 \\
 \geq & 2\varepsilon_0\|\bar{\varrho}_{0}\|_{L^2 } -\sqrt{C_4}\varepsilon_0^{3/2}
 \geq 2m_0\varepsilon_0 -\sqrt{C_4}\varepsilon_0^{3/2} \geq m_0\varepsilon_0,
 \end{aligned}      \end{equation*}
 Similar, we also have \begin{equation*}\begin{aligned}
 \|u_3^{\delta}(T^\delta)\|_{L^2 }
 \geq  2m_0\varepsilon_0 -\sqrt{C_4}\varepsilon_0^{3/2}\geq m_0\varepsilon_0,
 \end{aligned}      \end{equation*}
and
 \begin{equation*}\begin{aligned}
 \|(u_1^{\delta},u_2^{\delta})(T^\delta)\|_{L^2 }\geq
   2m_0\varepsilon_0 -\sqrt{C_4}\varepsilon_0^{3/2}
\geq m_0\varepsilon_0,
 \end{aligned}      \end{equation*}
 where $u^{\delta}_{i}(T^{\delta})$ denote the $i$-th component of
$ u^{\delta}(T^{\delta})$ for $i=1$, $2$, $3$.
This completes the proof of Theorem \ref{thm:0101new} by defining $\varepsilon :=m_0\varepsilon_0$.

\vspace{4mm} \noindent\textbf{Acknowledgements.} %The authors would
%like to thank the anonymous referee for invaluable suggestions,
%which improve the presentation of this paper.
The research of   Fei Jiang was supported by NSFC (Grant Nos. %11101044, 11271051,
11301083 and 11471134).% and the NSF of Fujian Province of China (Grant No. 2014J01011).

\renewcommand\refname{References}
\renewenvironment{thebibliography}[1]{%
\section*{\refname}
\list{{\arabic{enumi}}}{\def\makelabel##1{\hss{##1}}\topsep=0mm
\parsep=0mm
\partopsep=0mm\itemsep=0mm
\labelsep=1ex\itemindent=0mm
\settowidth\labelwidth{\small[#1]}%
\leftmargin\labelwidth \advance\leftmargin\labelsep
\advance\leftmargin -\itemindent
\usecounter{enumi}}\small
\def\newblock{\ }
\sloppy\clubpenalty4000\widowpenalty4000
\sfcode`\.=1000\relax}{\endlist}
\bibliographystyle{model1b-num-names}
%\bibliography{refs}

\begin{thebibliography}{24}
\expandafter\ifx\csname natexlab\endcsname\relax\def\natexlab#1{#1}\fi
\providecommand{\bibinfo}[2]{#2}
\ifx\xfnm\relax \def\xfnm[#1]{\unskip,\space#1}\fi
%Type = Article
\bibitem[{Bouya(2013)}]{BIIS413}
\bibinfo{author}{I.~Bouya}, \bibinfo{title}{{Instability of the forced
  magnetohydrodynamics system at small Reynolds number}},
  \bibinfo{journal}{SIAM J. Math. Anal.,} \bibinfo{volume}{45}
  (\bibinfo{year}{2013}) \bibinfo{pages}{307--323}.
%Type = Article
\bibitem[{Cho and Kim(2004)}]{CYKHOM1ufda}
\bibinfo{author}{Y.~Cho}, \bibinfo{author}{H.~Kim}, \bibinfo{title}{Unique
  solvability for the density-dependent Navier--Stokes equations},
  \bibinfo{journal}{Nonlinear Anal.} \bibinfo{volume}{59}
  (\bibinfo{year}{2004}) \bibinfo{pages}{465--489}.
%Type = Article
\bibitem[{Duan et~al.(2012)Duan, Jiang and Jiang}]{DRJFJS}
\bibinfo{author}{R.~Duan}, \bibinfo{author}{F.~Jiang},
  \bibinfo{author}{S.~Jiang}, \bibinfo{title}{{On the Rayleigh-Taylor
  instability for incompressible, inviscid magnetohydrodynamic flows}},
  \bibinfo{journal}{SIAM J. Appl. Math.} \bibinfo{volume}{71}
  (\bibinfo{year}{2012}) \bibinfo{pages}{1990--2013}.
%Type = Article
\bibitem[{Friedlander et~al.(2009)Friedlander, Nata$\mathrm{\check{s}}$a and
  Vicol}]{FSNPVVNC}
\bibinfo{author}{S.~Friedlander},
  \bibinfo{author}{P.~Nata$\mathrm{\check{s}}$a}, \bibinfo{author}{V.~Vicol},
  \bibinfo{title}{{Nonlinear instability for the critically dissipative
  quasi-geostrophic equation}}, \bibinfo{journal}{Commun. Math. Phys.}
  \bibinfo{volume}{292} (\bibinfo{year}{2009}) \bibinfo{pages}{797--810}.
%Type = Article
\bibitem[{Friedlander et~al.(1997)Friedlander, Strauss and Vishik}]{FSSWVMNA}
\bibinfo{author}{S.~Friedlander}, \bibinfo{author}{W.~Strauss},
  \bibinfo{author}{M.~Vishik}, \bibinfo{title}{{Nonlinear instability in an
  ideal fluid}}, \bibinfo{journal}{Ann. Inst. H. Poincare Anal. Non Lineaire}
  \bibinfo{volume}{14} (\bibinfo{year}{1997}) \bibinfo{pages}{187--209}.
%Type = Article
\bibitem[{Friedlander and Vishik(2003)}]{FrrVishikM}
\bibinfo{author}{S.~Friedlander}, \bibinfo{author}{M.~Vishik},
  \bibinfo{title}{{Nonlinear instability in two dimensional ideal fluids: the
  case of a dominant eigenvalue}}, \bibinfo{journal}{Comm. Math. Phys.}
  \bibinfo{volume}{243} (\bibinfo{year}{2003}) \bibinfo{pages}{261--273}.
%Type = Article
\bibitem[{G$\mathrm{\acute{e}}$rard-Varet(2005)}]{GDOS2}
\bibinfo{author}{D.~G$\mathrm{\acute{e}}$rard-Varet},
  \bibinfo{title}{{Oscillating solutions of incompressible magnetohydrodynamics
  and dynamo effect}}, \bibinfo{journal}{SIAM J. Math. Anal.}
  \bibinfo{volume}{37} (\bibinfo{year}{2005}) \bibinfo{pages}{815--840}.
%Type = Article
\bibitem[{Grenier(2000)}]{GRnierEmanu}
\bibinfo{author}{E.~Grenier}, \bibinfo{title}{On the nonlinear instability of
  euler and prandtl equations}, \bibinfo{journal}{Comm. Pure. Appl. Math.}
  \bibinfo{volume}{53} (\bibinfo{year}{2000}) \bibinfo{pages}{1067--1091}.
%Type = Article
\bibitem[{Guo and Han(2010)}]{GYHYQCQ}
\bibinfo{author}{Y.~Guo}, \bibinfo{author}{Y.Q. Han}, \bibinfo{title}{{Critical
  Rayleigh number in Rayleigh-B$\mm{\acute{e}}$nard convection}},
  \bibinfo{journal}{Quart. Appl. Math.} \bibinfo{volume}{68}
  (\bibinfo{year}{2010}) \bibinfo{pages}{149--160}.
%Type = Article
\bibitem[{Guo and Strauss(1995{\natexlab{a}})}]{GYSWIC}
\bibinfo{author}{Y.~Guo}, \bibinfo{author}{W.A. Strauss},
  \bibinfo{title}{Instability of periodic BGK equilibria},
  \bibinfo{journal}{Comm. Pure Appl. Math.} \bibinfo{volume}{48}
  (\bibinfo{year}{1995}{\natexlab{a}}) \bibinfo{pages}{861--894}.
%Type = Article
\bibitem[{Guo and Strauss(1995{\natexlab{b}})}]{GYSWICNonlinea}
\bibinfo{author}{Y.~Guo}, \bibinfo{author}{W.A. Strauss},
  \bibinfo{title}{Nonlinear instability of double-humped equilibria},
  \bibinfo{journal}{Ann. Inst. H. Poincar¡äe Anal. Non
  Lin$\mathrm{\acute{e}}$aire} \bibinfo{volume}{12}
  (\bibinfo{year}{1995}{\natexlab{b}}) \bibinfo{pages}{339--352}.
%Type = Article
\bibitem[{Guo and Tice(2011{\natexlab{a}})}]{GYTI1}
\bibinfo{author}{Y.~Guo}, \bibinfo{author}{I.~Tice},
  \bibinfo{title}{Compressible, inviscid Rayleigh--Taylor instability},
  \bibinfo{journal}{Indiana Univ. Math. J.} \bibinfo{volume}{60}
  (\bibinfo{year}{2011}{\natexlab{a}}) \bibinfo{pages}{677--712}.
%Type = Article
\bibitem[{Guo and Tice(2011{\natexlab{b}})}]{GYTI2}
\bibinfo{author}{Y.~Guo}, \bibinfo{author}{I.~Tice}, \bibinfo{title}{Linear
  Rayleigh--Taylor instability for viscous, compressible fluids},
  \bibinfo{journal}{SIAM J. Math. Anal.} \bibinfo{volume}{42}
  (\bibinfo{year}{2011}{\natexlab{b}}) \bibinfo{pages}{1688--1720}.
%Type = Article
\bibitem[{Hide(1955)}]{HRWP}
\bibinfo{author}{R.~Hide}, \bibinfo{title}{{Waves in a heavy, viscous,
  incompressible, electrically conducting fluid of variable density, in the
  presence of a magnetic field}}, \bibinfo{journal}{Proc. Roy. Soc. (London) A}
  \bibinfo{volume}{233} (\bibinfo{year}{1955}) \bibinfo{pages}{376--396}.
%Type = Article
\bibitem[{Hwang and Guo(2003)}]{HHJGY}
\bibinfo{author}{H.J. Hwang}, \bibinfo{author}{Y.~Guo}, \bibinfo{title}{{On the
  dynamical Rayleigh-Taylor instability}}, \bibinfo{journal}{Arch. Rational
  Mech. Anal.} \bibinfo{volume}{167} (\bibinfo{year}{2003})
  \bibinfo{pages}{235--253}.
%Type = Article
\bibitem[{Jiang and Jiang(2013)}]{JFJSOOADvace}
\bibinfo{author}{F.~Jiang}, \bibinfo{author}{S.~Jiang}, \bibinfo{title}{{ On
  instability and stability of three-dimensional gravity driven viscous flows
  in a bounded domain}}, \bibinfo{journal}{Adv. Math.} \bibinfo{volume}{264}
  (\bibinfo{year}{2013}) \bibinfo{pages}{831--863}.
%Type = Article
\bibitem[{Jiang et~al.(2013{\natexlab{a}})Jiang, Jiang and Ni}]{NJTSC2}
\bibinfo{author}{F.~Jiang}, \bibinfo{author}{S.~Jiang}, \bibinfo{author}{G.X.
  Ni}, \bibinfo{title}{{Nonlinear instability for nonhomogeneous incompressible
  viscous fluids}}, \bibinfo{journal}{Sci. China Math.} \bibinfo{volume}{56}
  (\bibinfo{year}{2013}{\natexlab{a}}) \bibinfo{pages}{665--686}.


%Type = Article
\bibitem[{Jiang et~al.(2013{\natexlab{b}})Jiang, Jiang and Wang}]{JFJSWWWN}
\bibinfo{author}{F.~Jiang}, \bibinfo{author}{S.~Jiang},
  \bibinfo{author}{W.~Wang}, \bibinfo{title}{{Nonlinear Rayleigh--Taylor
  instability in nonhomogeneous incompressible viscous magnetohydrodynamic
  fluids}}, \bibinfo{journal}{arXiv:1304.5636v1 [math.AP], 20 April 2013.}
  (\bibinfo{year}{2013}{\natexlab{b}}).
%Type = Article
\bibitem[{Jiang et~al.(2014)Jiang, Jiang and Wang}]{JFJSWWWOA}
\bibinfo{author}{F.~Jiang}, \bibinfo{author}{S.~Jiang},
  \bibinfo{author}{Y.~Wang}, \bibinfo{title}{{On the Rayleigh--Taylor
  instability for the incompressible viscous magnetohydrodynamic equations}},
  \bibinfo{journal}{Comm. Partial Differential Equation} \bibinfo{volume}{39}
  (\bibinfo{year}{2014}) \bibinfo{pages}{399--438}.
%Type = Article
\bibitem[{Kruskal and Schwarzschild(1954)}]{KMSMSP}
\bibinfo{author}{M.~Kruskal}, \bibinfo{author}{M.~Schwarzschild},
  \bibinfo{title}{Some instabilities of a completely ionized plasma},
  \bibinfo{journal}{Proc. Roy. Soc. (London) A} \bibinfo{volume}{233}
  (\bibinfo{year}{1954}) \bibinfo{pages}{348--360}.
%Type = Article
\bibitem[{Lin(2005)}]{LinzhuwuSS}
\bibinfo{author}{Z.~Lin}, \bibinfo{title}{Nonlinear instability of periodic BGK
  waves for vlasov-poisson system}, \bibinfo{journal}{Comm. Pure Appl. Math.}
  \bibinfo{volume}{LVIII} (\bibinfo{year}{2005}) \bibinfo{pages}{0505--0528}.
%Type = Article
\bibitem[{Pr$\mathrm{\ddot{u}}$ess and Simonett(2010)}]{PJSGOI5}
\bibinfo{author}{J.~Pr$\mathrm{\ddot{u}}$ess}, \bibinfo{author}{G.~Simonett},
  \bibinfo{title}{{On the Rayleigh-Taylor instability for the two-phase
  Navier-Stokes equations}}, \bibinfo{journal}{Indiana Univ. Math. J.}
  \bibinfo{volume}{59} (\bibinfo{year}{2010}) \bibinfo{pages}{1853--1871}.
%Type = Article
\bibitem[{Rayleigh(1883)}]{RLAP}
\bibinfo{author}{L.~Rayleigh}, \bibinfo{title}{{Analytic solutions of the
  Rayleigh equations for linear density profiles}}, \bibinfo{journal}{Proc.
  London. Math. Soc.} \bibinfo{volume}{14} (\bibinfo{year}{1883})
  \bibinfo{pages}{170--177}.
%Type = Article
\bibitem[{Wang and Tice(2012)}]{wang2011viscous}
\bibinfo{author}{Y.~Wang}, \bibinfo{author}{I.~Tice}, \bibinfo{title}{The
  viscous surface-internal wave problem: nonlinear Rayleigh--Taylor
  instability}, \bibinfo{journal}{Comm. in Partial Differential Equations}
  \bibinfo{volume}{37} (\bibinfo{year}{2012}) \bibinfo{pages}{1967--2028}.

\end{thebibliography}

\end{document}